\documentclass[a4paper,12pt]{article}
\usepackage[hideboxes]{boxedeps}

\usepackage{latexsym}
\newcommand{\pic}[2]{\BoxedEPSF{#1 scaled #2}}

\def\CC{{\cal C}}
\def\CS{{\cal S}}

\def\Xor{\pic{Xor.ART} {500}}
\def\Yor{\pic{Yor.ART} {500}}
\def\Ior{\pic{Ior.ART} {500}}
\def\Sigmaior{\pic{Sigmaior.ART} {500}}

\def\Idor{\pic{Idor.ART} {500}}
\def\Rcurlor{\pic{Rcurlor.ART} {500}}

\def\XC{\pic{XC.ART} {500}}
\def\XCYC{\pic{XCYC.ART} {500}}
\def\XYC{\pic{XYC.ART} {500}}
\def\TC{\pic{TC.ART} {500}}
\def\PhiX{\pic{PhiX.ART} {500}}
\def\PhiOne{\pic{PhiOne.ART} {500}}

\def\PsiX{\pic{PsiX.ART} {500}}

\def\HopfXY{\pic{HopfXY.ART} {500}}
\def\HopfXHopfY{\pic{HopfXHopfY.ART} {500}}
\def\Threeone{\pic{Threeone.ART}{500}}
\def\Twotwo{\pic{Twotwo.ART}{500}}

\def\Hopf{\pic{Hopf.ART} {500}}
\def\Hopfcurl{\pic{Hopfcurl.ART} {500}}
\def\Hopfcurlpar{\pic{Hopfcurlpar.ART} {500}}

\newcommand{\hopf}[2]{\langle #1,#2 \rangle}
\newcommand{\eval}[1]{\left\langle #1 \right\rangle}
\newcommand{\trans}[1]{{#1}^{\vee}}

\newcommand{\emptyd}{\emptyset}

\newcommand{\mylabel}[1]{\label{#1}}
\newcommand{\mysln}{sl(N)}

\newtheorem{theorem}{Theorem}[section]
\newtheorem{corollary}[theorem]{Corollary}
\newtheorem{lemma}[theorem]{Lemma}

\newenvironment{remark}{\par\smallskip%
\noindent\textbf{Remark.}\  }%
{\par\smallskip}

{\par\smallskip}

\newenvironment{definition}{\par\smallskip%
\noindent\textbf{Definition.}\  }%
{\par\smallskip}

\newenvironment{proof}[1][{}]{\par\smallskip%
\noindent\textit{Proof #1: }\  }
{\hfill$\Box$\par\smallskip}

\begin{document}
\begin{center}
      {\Large\bf The Homfly polynomial of the decorated Hopf link}

\bigskip
HUGH R. MORTON and SASCHA G. LUKAC
\footnote{Funded by DAAD under their HSP-III programme}

\medskip
 {\it Department of Mathematical Sciences, University of Liverpool,\\
     Peach St, Liverpool, L69 7ZL, England. }\\
{\tt morton@liv.ac.uk, lukac@liv.ac.uk}
\end{center}
\vskip 1cm

\begin{abstract}
 The main goal is to find the Homfly polynomial of a link formed by decorating
each component of the Hopf link with the closure of a directly oriented
tangle. Such decorations are spanned in the Homfly skein of the annulus by
elements $Q_\lambda$, depending on partitions $\lambda$. We show how
the 2-variable Homfly invariant $\hopf{\lambda}{\mu}$ of the Hopf link arising
from decorations $Q_\lambda$ and $Q_\mu$ can be found from the Schur symmetric
function $s_\mu$ of an explicit power series depending on $\lambda$. We show
also that the quantum  invariant of the Hopf link coloured by irreducible
$sl(N)_q$ modules $V_\lambda$ and $V_\mu$, which is  a 1-variable
specialisation of $\hopf{\lambda}{\mu}$, can be expressed in terms of an
$N\times N$ minor of the Vandermonde matrix $(q^{ij})$.

\phantom{some stuff to separate the keywords from the abstract}

\end{abstract}

\section*{Introduction.}

The roots of this paper lie in the work of Morton and Strickland,
\cite{MorStr}, where the central role of the Hopf link in studying the
invariants of satellite knots became clear in the context of the coloured
Jones invariants. The Hopf link invariants also play a crucial part in the
construction of 3-manifold invariants based on the Jones polynomial. As
noted in \cite{MorStr}, the behaviour of the Hopf link invariants for
quantum groups such as  $\mysln_q$ was anticipated to be at the heart of
constructions of corresponding 3-manifold invariants for other values of
$N>2$, as borne out  by subsequent work such as that of Kohno and Takata
\cite{Kohno}. 

It is of considerable interest to find a
general formula for the Homfly polynomial of the Hopf link in which each
component is decorated by, for example, a closed braid. This will be a Laurent
polynomial in the parameters
$v$ and $z$, as used in
\cite{AistonMorton}, whose evaluation at the `$\mysln$ specialisation', with
$v=s^{-N}$ and $z=s-s^{-1}$, gives a Laurent polynomial in $q=s^2$ which can
be viewed in terms of
$\mysln_q$ invariants of the Hopf link.

In this paper we give a formula in theorem \ref{Frobenius}  to determine the
Homfly polynomial
$\hopf{\lambda}{\mu}$ of the Hopf link with components decorated by the
closures
$Q_\lambda, Q_\mu$ of the Gyoja-Aiston idempotents, \cite{Gyoja, Aiston},
corresponding to any two partitions $\lambda, \mu$. Invariants determined using these
decorations behave simply under change of framing.
The set $\{Q_\lambda\}$ spans a large subspace in the Homfly skein of the
annulus, including all elements represented by closed braids or more general
closed tangles without reverse strings.  Our calculations therefore
will find in principle the Homfly polynomial of satellites of the Hopf link
with any such decorations and any framing.

   As a key step in establishing the formula we present an
attractive expression in theorem
\ref{LambdaMu}  for the
$\mysln$ specialisation of
$\hopf{\lambda}{\mu}$ as a Laurent polynomial in the single variable $q$ in
terms of an $N\times N$ minor of the Vandermonde matrix $(q^{ij})$, when
$\lambda$ and $\mu$ each have at least $N$ parts. This Laurent polynomial in
$q$ also determines, up to an explicit power of $q$, the $\mysln_q$ invariant
of the Hopf link when its components are coloured by the irreducible
$\mysln_q$ modules $V_\lambda$ and $V_\mu$.

Much of our work was originally inspired by a similar expression for the
further specialisation of this determinant to $q=\exp(2\pi i/r)$ used by Kohno
and Takata \cite{Kohno}. Our proofs for generic
$v$ and
$s$ do not draw on their work and hence give an alternative skein theoretic 
way to interpret their formulae.

To describe and analyse the invariants under consideration we apply Homfly
skein theory to the skein of the annulus and the Hecke algebras, coupled with
methods from the classical theory of symmetric functions. We start with a
brief review of this material before applying it to the Hopf link.

Much of the material in this paper appeared originally in the thesis of Lukac,
\cite{Lukac}.

\section{The skein models.}
The account here largely follows those of \cite{AistonMorton}, \cite{Lukac}
and \cite{Murphy}.
 The  framed Homfly skein relations, in
their simplest form, are
\begin{eqnarray*}
\Xor\  -\ \Yor& =&(s-s^{-1})\ \Ior\\
\mbox{and }\  \Rcurlor &=&v^{-1}\ \Idor.
\end{eqnarray*}
The coefficient ring can be taken
as $\Lambda={\bf Z}[v^{\pm1}, s^{\pm1}]$ with
monomials in $\{s^k-s^{-k}:k\ge0\}$ admitted as denominators.
The framed Homfly skein ${\cal S}(F)$ of a planar
surface $F$, with some designated input and
output boundary points, is defined to be $\Lambda$-linear
combinations of oriented tangles in $F$, modulo
these two local relations, and Reidemeister moves
II and III. 

The empty tangle is admitted when $F$ has no boundary points. We include the
local relation which allows the removal of a null-homotopic closed curve
without crossings, at the expense of multiplication by the scalar
$\delta=\displaystyle\frac{v^{-1}-v}{s-s^{-1}}$. This is in fact a consequence
of the main relations, except when removing the curve leaves only the empty
diagram.

\subsection{The plane}
When $F={\bf R}^2$ every element can be represented uniquely as a scalar
multiple of the empty diagram. For a diagram $D$ the resulting scalar
$\chi(D)\in \Lambda$ is the framed Homfly polynomial of $D$. Taking
$\chi^u(D)=v^{wr(D)}\chi(D)$, where $wr(D)$ is the writhe of the diagram,
gives a scalar which is invariant under all Reidemeister moves. It is the
Homfly polynomial of $D$, as a function of $v$ and $s$, defined by the
local relation
\[v^{-1}\Xor\ -\ v\Yor=(s-s^{-1})\ \Ior,\]and normalised to have the value 1
on the empty diagram. Then $\chi^u(D)=\delta P(D)$, except when $D$ is empty,
where $P$ is the more traditional Homfly polynomial defined as the ambient
isotopy invariant which satisfies the local
relation above and takes the value $1$ on the unknot .

\subsection{The Hecke algebras}

Write $R_n^n$ for the skein ${\cal S}(F)$ of
$n$-tangles,  where $F$
is a rectangle with $n$ inputs at the bottom and
$n$ outputs at the top. Composing $n$-tangles
induces a product which makes $R_n^n$ into
an algebra. It has a linear basis of $n!$
elements, and is isomorphic to the Hecke algebra
$H_n(z)$, with coefficients extended to the ring $\Lambda$. This algebra has a
presentation generated by the  elementary
braids
\[\sigma_i\ =\ \Sigmaior\] subject to 
 the
braid relations
\begin{eqnarray*}
\sigma_i\sigma_j&=&\sigma_j\sigma_i,
\quad |i-j|>1,\\
\sigma_i\sigma_{i+1}\sigma_i&=&\sigma_{i+1}\sigma_i\sigma_{i+1},
\end{eqnarray*} and the quadratic relations
$\sigma_i^2=z\sigma_i+1$, with $z=s-s^{-1}$, giving the
alternative form $(\sigma_i-s)(\sigma_i+s^{-1})=0$.

A simple adjustment of the skein relations, as in \cite{AistonMorton},
allows for a skein model of $H_n$ whose parameters can be readily
adapted to match any of the different appearances of the algebra,
\cite{Murphy}.

\subsection{The  annulus.}

The Homfly
skein of the annulus, ${\cal C}$, as discussed in
\cite{Turkey} and originally in \cite{Turaev},  is defined to be linear
combinations of diagrams in the annulus, modulo the framed Homfly skein
relations. An element $X\in\CC$ will be indicated on a diagram as
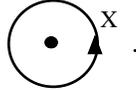
\[\XC\ .\]
The skein ${\cal C}$ has a product induced by
placing one annulus outside another,  under which $\cal C$ becomes a
commutative algebra;
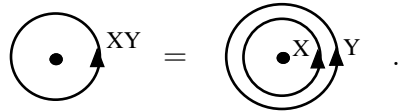
\[\XYC\ =\ \XCYC\ .\]

There is an evaluation map $\eval{\ }:\CC\to\Lambda$, induced by the inclusion
of a standard annulus in the plane, in which $\eval{X}$ is defined as the
framed Homfly polynomial of $X$ when regarded as a diagram in the plane. Since
the framed Homfly polynomial with our normalisation is multiplicative on split
diagrams, the evaluation map $\eval{\ }$ is a
ring homomorphism.

\subsection{Interrelations}

The best known relation of $\CC$ with the Hecke algebra $H_n$ is the
closure map
$R_n^n\to {\cal C}$, induced by taking a tangle $T$ to its
closure $\hat {T}$ in the annulus, defined by \[\hat {T}\ =\ \TC\ .\]
This is a linear map, whose image we call ${\cal C}_n$.

Elements of the skein $\CC$ can be used to decorate closed curves in a framed
diagram in any $F$ to produce new skein elements of $\CS(F)$.  We shall
restrict the decorations in our work to elements drawn from $\CC_n$ for some
$n$, and write $\displaystyle \CC^+=\cup_{n\ge0}\CC_n$.

A  simple skein theory construction
determines a natural linear map $\varphi:{\cal
C}\to{\cal C}$, induced by taking a diagram $X$
in the annulus and linking it once with a simple
loop to get
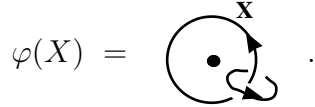
\[\varphi(X)\ =\ \PhiX\ .\]

This construction can be interpreted as an example of the decoration technique
above, applied to one loop of the diagram \[\PhiOne.\]
 We can extend this  
by decorating the other loop with any
$Y\in\CC$ to give a linear map  $\varphi(Y):\CC\to\CC$. Then
\[\varphi(Y)(X)\ =\ \HopfXY\ \in\CC.\]
 The  map
$\varphi$ is then $\varphi(c_1)$ where $c_1\in\CC$ is 
represented by a single closed curve following the oriented core of the
annulus.

Evaluation of $\varphi(Y)(X)$  gives a scalar which we
 denote  by $\hopf{X}{Y}$, and which forms the main focus of this
paper.

The map $\psi_n:\CC\to H_n$ induced by decorating the closed curve in the
diagram below  by an element $X\in\CC$ is readily seen to be a ring
homomorphism, whose image
\[\psi_n(X)\ =\ \PsiX\ \in H_n\]
lies in the centre of $H_n$.

There is a  set of quasi-idempotent elements,
$e_\lambda\in H_n$, one for each partition $\lambda$ of
$n$, known as the Gyoja-Aiston idempotents.  They were originally described
algebraically by Gyoja \cite{Gyoja}, while skein pictures of
these based on the Young diagram for $\lambda$
can be found in \cite{Aiston} or \cite{AistonMorton}. The skein theory version
displays a nice `internal stability' under multiplication which shows readily
that $Ze_\lambda$ is a scalar multiple of $e_\lambda$ for any {\em central}
element
$Z\in H_n$. We can then define a ring homomorphism $t_\lambda:\CC\to
\Lambda$ by the formula $\psi_n(Y)e_\lambda=t_\lambda(Y)e_\lambda$, using the
homomorphism $\psi_n$ from $\CC$ to the centre of $H_n$.

Define $Q_\lambda\in \CC$ to be the closure of the genuine idempotent
$\frac{1}{\alpha_\lambda}e_\lambda$, where $e_\lambda^2=\alpha_\lambda
e_\lambda$, as in \cite{AistonMorton} or \cite{Lukac}. Then clearly
$Q_\lambda$ is an eigenvector of $\varphi(Y)$ with eigenvalue $t_\lambda(Y)$
for every $Y\in\CC$. We make extensive use of these elements of $\CC$,
especially where $\lambda$ is a single column or row. We write $c_i$ for
$Q_\lambda$ where $\lambda$ is the column with $i$ cells, and $d_j$ where
$\lambda$ is the row with $j$ cells, and take $c_0=d_0=1\in\CC$ represented
by the empty diagram.

\subsection{The Hopf link}

We consider the Hopf link with linking number 1 as shown here.

  \begin{center}
    \Hopf
     
  \end{center}

Let $X$ and $Y$ be any elements of the skein $\CC$ of the annulus.
We write $\hopf{X}{Y}$ for
 the framed Homfly polynomial of the
Hopf link with decorations $X$ and $Y$ on its components,
giving\[\hopf{X}{Y}=\chi\left(\HopfXHopfY\right)\ =\ \eval\HopfXY\
=\eval{\varphi(Y)(X)}.\]
 Clearly $\hopf{Y}{X}=\hopf{X}{Y}$. Our aim is to
determine
$\hopf{X}{Y}$ for any $X,Y\in\CC^+$. Since $\CC^+$ is spanned by
$\{Q_\lambda\}$ it is enough to find $\hopf{Q_\lambda}{Q_\mu}$ for all
$\lambda,\mu$. We make use of the homomorphism $t_\lambda$ of the previous
section.

\begin{lemma}\label{C-etalamhom}
We can write \[t_\lambda(Y)=\hopf{Q_\lambda}{Y}/\eval{Q_\lambda}.\] 
\end{lemma}
\begin{proof}
Suppose that $\lambda$ has $n$ cells. Then
$\psi_n(Y)e_\lambda=t_\lambda(Y)e_\lambda$. Apply the closure map
$\hat{\ }:H_n\to
\CC$ to both sides. The left-hand side becomes $\varphi(Y)(\hat{e}_\lambda)$.
Now apply the evaluation map
$\eval{\ }$ to both sides, using the fact that
$\hat{e}_\lambda=\alpha_\lambda Q_\lambda$, to get
\[\hopf{\alpha_\lambda Q_\lambda}{Y}=t_\lambda(Y)\eval{\alpha_\lambda
Q_\lambda}.\]
This gives the expression for $t_\lambda(Y)$, since $\alpha_\lambda\ne0$.
\end{proof}

\begin{corollary}\label{L-hopf2prod}
For any elements $X$ and $Y$ of $\CC$ and any Young diagram $\lambda$
we have
\[
\hopf{Q_{\lambda}}{X}\hopf{Q_{\lambda}}{Y}=
\eval{Q_{\lambda}}\hopf{Q_{\lambda}}{XY}.
\]
\end{corollary}
\begin{proof}
Use the formula of lemma \ref{C-etalamhom} and the fact that $t_\lambda$ is a
homomorphism.
\end{proof}

Since any $Y\in\CC^+$ can be written as a polynomial in the skein elements
$\{c_i\}$ we only need to know the values of
$\hopf{Q_\lambda}{c_i}$ for  integers $i\geq 0$ in order to
compute
$\hopf{Q_\lambda}{Y}$ for all $Y\in\CC^+$. In particular we can then compute
all invariants $\hopf{Q_\lambda}{Q_\mu}$.
  Hence, it is useful to define a
formal power series
\begin{equation}
  E_{\lambda}(t)=
\sum_{i\geq 0}t_\lambda(c_i)t^i=\frac{1}{\eval{Q_\lambda}}\sum_{i\geq
0}\hopf{Q_\lambda}{c_i}t^i \label{Elambda}
\end{equation}
for any Young diagram $\lambda$, which we study further in the next section.

\section{Symmetric functions and the skein of the annulus.}

In this section we recall some explicit results about elements in the
Hecke algebras and their closure in $\CC$, and their interpretation in
the context of symmetric functions, following the methods of
Macdonald \cite{Mac}.  

\subsection{Formal Schur functions}

We  follow Macdonald in defining  the Schur functions of a formal power
series $E(t)=1+\sum_{i=1}^\infty  e_it^i$ with coefficients $e_i$ in a
commutative ring. When $E(t)$ is a polynomial with a formal
factorisation as $E(t)=\prod_{j=1}^N (1+x_jt)$ then any partition $\lambda$ of
$n$ determines the classical Schur symmetric function
$s_\lambda(x_1,\ldots,x_N)$ in terms of determinants whose
entries are powers of $\{x_j\}$, given explicitly by equation
(\ref{Schurclass}). This can be written as a polynomial in the elementary
symmetric functions
$\{e_i\}$ of
$x_1,\ldots,x_N$, which are the coefficients of $E(t)$, by means of the
Jacobi-Trudy formula. So long as $N\ge n$ this polynomial in $\{e_i\}$ is
independent of $N$. Macdonald  defines the Schur function
$s_\lambda(E(t))$ for a formal power series to be this polynomial in the
coefficients $\{e_i\}$, described explicitly in section \ref{JacobiTrudy}.

 The coefficient $e_i$ itself is the Schur function
$s_\lambda(E(t))$ for the partition $\lambda =1^i=(1,\ldots,1)$, represented
by the Young diagram consisting of a single column with $i$ cells. Equally,
when
$\lambda$ is a single row with $j$ cells then $s_\lambda(E(t))=h_j$, where
$h_j$ is the complete symmetric function of degree $j$ in $x_1,\ldots,x_N$ in
the case of the polynomial $E(t)$ above. In any case, the power series
$H(t)=1+\sum h_jt^j$ satisfies the power series equation $E(-t)H(t)=1$.

Since $s_\lambda(E(t))$ is a polynomial in $\{e_i\}$ depending only on
$\lambda$ it is clear that if $R$ and $S$ are any commutative rings and we
apply a ring homomorphism
$\rho:R\to S$ to the coefficients $\{e_i\}\in R$ to get the power series
$\rho(E(t))$ with coefficients $\{\rho(e_i)\}\in S$ then
\[s_\lambda(\rho(E(t)))=\rho (s_\lambda(E(t))),\] for every $\lambda$.

We make use later of the homogeneity of the Schur function $s_\lambda$ in the
form
\begin{equation}s_\lambda(E(\alpha t))=\alpha^{|\lambda|}s_\lambda(E(t)),
\label{Homogeneous}
\end{equation}
where $\lambda$ has $|\lambda|$ cells.

\subsection{Applications in the skein of the annulus}

The elements $Q_\lambda$ have a very nice interpretation  as Schur
functions in $\CC$.  Lukac showed in \cite{Lukac} that $Q_\lambda$ can be identified with the
Schur function
\[s_\lambda(\sum_{i\ge0} c_it^i)\] of the series  $\sum c_it^i$ whose
coefficients $c_i\in\CC$ are the closures $Q_\lambda$ of the single column
idempotents with $\lambda =1^i=(1,\ldots,1)$, under the convention that
$c_0=1\in\CC$, represented by the empty diagram.  This has been known in
principle for some time,  but the proof in
\cite{Aiston} was fairly circuitous and had to refer at one point to  quantum
group work. Kawagoe
\cite{Kawagoe} used these Schur functions and proved that they were
eigenvectors of
$\varphi$, but his proof too was relatively complicated. Lukac gave a much
more satisfactory skein proof, using a further skein-based algebra, coupled
with the knowledge that the eigenvalues of  $\varphi|\CC_n$ are distinct and
the elements
$Q_\lambda$ can essentially be characterised as eigenvectors of $\varphi$.

The fact that $Q_\lambda$ is also an eigenvector of every $\varphi(Y)$, which
is not at all clear at first sight for the Schur function $s_\lambda(\sum
c_it^i)$, is the feature that allows us to complete our analysis here, using
the consequence that $t_\lambda$ is a homomorphism. 

For the sake of brevity of notation we write $\lambda$ in place of
$Q_\lambda$ as an element of $\CC$. As Schur functions of a
power series these elements  multiply according to the Littlewood-Richardson
rules for Young diagrams. We make use of this later when we 
calculate the product in $\CC$ of the elements $c_i$ and $d_j$ corresponding to
a single row and column respectively.

From the series $E_\lambda(t)$, defined in (\ref{Elambda}) above,  we
can find
$\hopf{\lambda}{\mu}$ by calculating its Schur function $s_\mu(E_\lambda(t))$
and using the following lemma.
\begin{lemma}\mylabel{L-muschur}
We have
\[
  s_{\mu}(E_{\lambda}(t))=\frac{1}{\eval{\lambda}}\hopf{\lambda}{\mu}
\]
for any Young diagrams $\lambda$ and $\mu$.
\end{lemma}

\begin{proof} Applying the homomorphism $t_\lambda$ to the power series
$\sum c_it^i$ gives $E_\lambda(t)=t_\lambda(\sum c_it^i)$. Then \[
s_\mu(E_\lambda(t))=t_\lambda(s_\mu(\sum
c_it^i))=t_\lambda(\mu)=\frac{1}{\eval{\lambda}}\hopf{\lambda}{\mu}.
\]
\end{proof}

We also define 
\[
  H_{\lambda}(t)=\frac{1}{\eval{\lambda}}\sum_{j\geq
0}\hopf{\lambda}{d_j}t^j=t_\lambda(\sum d_jt^j)
\]
for any Young diagram $\lambda$, taking $d_0=1$.

\begin{lemma}\label{L-hlamelam}
We have
\[
  E_{\lambda}(t) H_{\lambda}(-t)=1
\]
for any Young diagram $\lambda$.
\end{lemma}
\begin{proof} Apply the homomorphism $t_\lambda$ to the two power 
series with coefficients in $\CC$ in the equation  \[
  \left(\sum_{i\geq 0} c_i t^i\right)\left(\sum_{j\geq 0} d_j (-t)^j\right)=1
\] proved by direct skein theory in \cite{Aiston}. The equation follows
alternatively from the fact that the elements $d_j$ in $\CC$ are the Schur
functions corresponding to single row diagrams.
\end{proof}

\section{The Hopf link decorated with columns and rows}

We now compute the series $E_{c_k}(t)$ for any integer $k\geq 0$.
To do this, we start with a surprisingly simple formula for
$\hopf{c_i}{d_j}$.
\begin{lemma}\label{L-hopfckdj}
We have
\[
  \hopf{c_i}{d_j}=\eval{c_i}\eval{d_j}
                  \frac{v^{-1}(s^{2j}-s^{2(j-i)}+s^{-2i})-v}{v^{-1}-v}
\]
for any integers $i\geq 0$ and $j\geq 0$.
\end{lemma}
\begin{proof}
We shall prove the lemma by expressing the Homfly polynomial
of a certain decorated link
in two different ways and comparing the results.
The link in question is the 2-parallel of the unknot with framing 1
as shown below,
decorated with ${c_i}$ on one component and ${d_j}$ on the other component.
We denote its framed Homfly polynomial by $R$.

  \begin{center}
    $R\ =\ \chi\left(\Hopfcurlpar\right)\ =\ \chi\left(\Hopfcurl\right)$
      \end{center}

The positive curl on $n$ strings belongs to the centre of the Hecke algebra
$H_n$, and so its product with any
quasi-idempotent
$e_{\lambda}$,
$|\lambda|=n$,
is a scalar multiple  of $e_{\lambda}$. The {\em content}, $cn(x)$, of the
cell $x$ in row $i$ and column $j$ of a Young diagram $\lambda$ is defined to
be $cn(x)=j-i$. The scalar
$f(\lambda)$ was calculated using skein theory in theorem 17 of
\cite{AistonMorton} as
\[f(\lambda)=v^{-|\lambda|}s^{n_{\lambda}}\]
where $n_{\lambda}$ is twice the sum of the contents of all cells
of $\lambda$.  
By removing the
 curls on the two components of the 2-parallel  we get
\begin{eqnarray}\label{E-rfck}
  R=f(c_i)f(d_j)\hopf{c_i}{d_j}.
\end{eqnarray}

The other way to calculate $R$ is to regard it as the Homfly 
polynomial of the unknot with framing 1 decorated by the product
of ${c_i}$ and ${d_j}$, as elements of $\CC$.
Since the elements $Q_{\lambda}$ are Schur functions, by Lukac
\cite{Lukac}, they multiply according to the Littlewood-Richardson rules. We
can then write
${c_i}{d_j}={\mu_{i,j+1}}+{\mu_{i+1,j}}$ where
$\mu_{a,b}$ is the simple hook Young diagram with 
$a$ cells in the first column and $b$ cells in the first row.
Hence
\begin{eqnarray}\label{E-rfmu}
  R=f(\mu_{i,j+1})\eval{\mu_{i,j+1}}+f(\mu_{i+1,j})\eval{\mu_{i+1,j}}.
\end{eqnarray}

From the  formula for $f(\lambda)$ above we get
\[f(c_i)=v^{-i}s^{-i(i-1)}, f(d_j)=v^{-j}s^{j(j-1)},
f(\mu_{a,b})=v^{-(a+b-1)}s^{b(b-1)-a(a-1)},\] and hence
\[
  f(c_{i+1})=v^{-1}s^{-2i}f(c_i),
  \ f(d_{j+1})=v^{-1}s^{2j}f(d_j),\ 
  f(\mu_{i,j})=vf(c_i)f(d_j).
\]

The following  relations were shown skein theoretically in
\cite{Aiston}. 
\begin{eqnarray*}
  \eval{c_{i+1}}&=&\frac{v^{-1}s^{-i}-vs^i}{s^{i+1}-s^{-i-1}}\eval{c_i},\\
  \eval{d_{j+1}}&=&\frac{v^{-1}s^j-vs^{-j}}{s^{j+1}-s^{-j-1}}\eval{d_j},\\
  \eval{\mu_{i,j}}&=&\frac{(s^j-s^{-j})(s^i-s^{-i})}
                          {(v^{-1}-v)(s^{i+j-1}-s^{-i-j+1})}
                     \eval{c_i}\eval{d_j}.
\end{eqnarray*}

We can in fact  deduce the product formula (\ref{Product}) for
$\sum\eval{c_i}t^i$ using only the first of these, and hence the closed
formula (\ref{Schurprod}) for
$\eval{\lambda}$.

 Using these relations equation (\ref{E-rfmu}) gives
\begin{eqnarray}
R &=& f(\mu_{i,j+1})\eval{\mu_{i,j+1}}+f(\mu_{i+1,j})\eval{\mu_{i+1,j}}\nonumber\\
  &=& vf(c_i)v^{-1}s^{2j}f(d_j)
       \frac{(s^i-s^{-i})(v^{-1}s^j-vs^{-j})}
            {(v^{-1}-v)(s^{i+j}-s^{-i-j})}
                     \eval{c_i}\eval{d_j}\nonumber\\
  & & +{}s^{-2i}f(c_i)f(d_j)\frac{(s^j-s^{-j})(v^{-1}s^{-i}-vs^i)}
                  {(v^{-1}-v)(s^{i+j}-s^{-i-j})}\eval{c_i}\eval{d_j}\nonumber\\
  &=&\frac{v^{-1}(s^{2j}-s^{2(j-i)}+s^{-2i})-v}{v^{-1}-v}
     f(c_i)f(d_j)\eval{c_i}\eval{d_j}.\label{E-rfck2}
\end{eqnarray}
 Equations
(\ref{E-rfck}) and (\ref{E-rfck2}) then give
\[
  \hopf{c_i}{d_j}=\eval{c_i}\eval{d_j}
          \frac{v^{-1}(s^{2j}-s^{2(j-i)}+s^{-2i})-v}{v^{-1}-v}
\] since $f(c_i)$ and $f(d_j)$ are non-zero.
\end{proof}

Using our notation above, we have
\begin{corollary}\label{C-hckprod}
\[
  H_{c_k}(t)=\frac{1-v^{-1}s^{-2k+1}t}{1-v^{-1}st}
             H_{\emptyd}(t)
\]
for any integer $k\geq 0$.
\end{corollary}
\begin{proof}
Since
$H_{\emptyd}(t)=t_\emptyd(H(t))=\sum_{j\geq0}\eval{d_j}t^j$ we must  show
that
\[
  (1-v^{-1}st)\frac{1}{\eval{c_k}}\sum_{j\geq 0}\hopf{c_k}{d_j}t^j
  =(1-v^{-1}s^{-2k+1}t)\sum_{j\geq 0}\eval{d_j}t^j.
\]
The constant terms of the power series in the above equation are
equal to 1. To show that the coefficients of $t^j$ on each side
agree it is enough
to show that
\begin{eqnarray}\label{E-hopfcoeff}
  \frac{1}{\eval{c_k}}\hopf{c_k}{d_j}-v^{-1}s\frac{1}{\eval{c_k}}
    \hopf{c_k}{d_{j-1}}=\eval{d_j}-v^{-1}s^{-2k+1}\eval{d_{j-1}}.
\end{eqnarray}
By lemma \ref{L-hopfckdj}  the left hand side of
equation (\ref{E-hopfcoeff}) is
\begin{eqnarray*}
& &
  \eval{d_j}\frac{v^{-1}(s^{2j}-s^{2(j-k)}+s^{-2k})-v}{v^{-1}-v}\\
& &
  -v^{-1}s\eval{d_{j-1}}\frac{v^{-1}(s^{2(j-1)}-s^{2(j-1-k)}+s^{-2k})-v}
          {v^{-1}-v}.
\end{eqnarray*}
Because
\[
  \eval{d_j}=\eval{d_{j-1}}\frac{v^{-1}s^{j-1}-vs^{-j+1}}{s^j-s^{-j}}
\]
this can be rewritten as
\begin{eqnarray}
  & &\left(
     \frac{(v^{-1}s^{j-1}-vs^{-j+1})(v^{-1}(s^{2j}-s^{2(j-k)}+s^{-2k})-v)}
       {(s^j-s^{-j})(v^{-1}-v)}\right.\nonumber\\
  & &\ \left. \ -v^{-1}s\frac{
              v^{-1}(s^{2(j-1)}-s^{2(j-1-k)}+s^{-2k})-v}{v^{-1}-v}
  \right)\eval{d_{j-1}}.\label{E-hopfcoeff2}
\end{eqnarray}  
The right hand side of equation (\ref{E-hopfcoeff}) is 
\begin{eqnarray}\label{E-hopfcoeff3}
  \left(
        \frac{v^{-1}s^{j-1}-vs^{-j+1}}{s^{j}-s^{-j}} - v^{-1}s^{-2k+1}
  \right)\eval{d_{j-1}}.
\end{eqnarray}
It is straightforward to confirm that the expressions (\ref{E-hopfcoeff2})
and (\ref{E-hopfcoeff3}) are equal, and thus equation (\ref{E-hopfcoeff})
follows.
\end{proof}
An immediate consequence of corollary \ref{C-hckprod}
and lemma \ref{L-hlamelam} is
\begin{corollary}\label{C-ckpower}
\[
  E_{c_k}(t)=\frac{1+v^{-1}st}{1+v^{-1}s^{-2k+1}t}
             E_{\emptyd}(t)
\]
for any integer $k\geq 0$.
\end{corollary}

\section{The Hopf link decorated with any Young diagrams}

We now in principle have a means of finding $\hopf{c_k}{\lambda}$ by
calculating the Schur function $s_\lambda(E_{c_k}(t))$. This in turn gives the
coefficients for the series $E_\lambda(t)$, whose Schur function $s_\mu$
finally gives $\hopf{\lambda}{\mu}$. The aim of this section is to give a
simple formula for the series
$E_\lambda(t)=\sum\hopf{\lambda}{c_i}t^i/\eval{\lambda}$ as a product of the
known series
$E_\emptyd(t)=\sum\eval{c_i}t^i$ and a rational function of $t$, in theorem
\ref{Frobenius}. To do this we specialise the coefficients so that
$E_{c_k}(t)$ becomes a polynomial, and we can apply the 
classical determinantal formulae for Schur functions. By comparing the
coefficients  in $E_\lambda(t)$ with those of our proposed product formula,  
 and showing that they agree under sufficiently many of the specialisations we
are able to deduce that they are identical.

The determinantal formulae which appear give us a very attractive expression
for the specialised versions of the invariants $\hopf{\lambda}{\mu}$, which we
describe in theorem \ref{LambdaMu}.

\subsection{The $\mysln$ substitution}

We use the substitution $v=s^{-N}$ to define a ring homomorphism from our
coefficient ring $\Lambda$ to the ring of Laurent polynomials in $s$ with
denominators $s^r-s^{-r}$. Denote the image of  $\hopf{X}{Y}$ by
$\hopf{X}{Y}_N$, and the image of $\eval{X}$ by $\eval{X}_N$.

The element $t_\lambda(Y)$ for $Y\in\CC$ is an element of $\Lambda$.
Applying the homomorphism to the equation
$\eval{\lambda}t_\lambda(Y)=\hopf{\lambda}{Y}$ from lemma \ref{C-etalamhom}
shows that if
$\eval{\lambda}_N=0$ then $\hopf{\lambda}{Y}_N=0$ also, for any $Y$.

We  make use of some results of Macdonald to help in our evaluations.
In  exercise I.3.3 of \cite{Mac} Macdonald observes that when a power series
can be expressed as
\[E(t) =\prod_{i=0}^\infty \frac{1+aq^it}{1+bq^it}\] then its Schur function
$s_\lambda$ is given in terms of the cells $x\in\lambda$ by the formula 
\begin{equation}
s_\lambda=q^{n(\lambda)}\prod_{x\in
\lambda}\frac{a-bq^{cn(x)}}{1-q^{hl(x)}},\label{Schurmac}
\end{equation}
where $cn(x)$ is the content of $x$, determined by its position, $hl(x)$ is its
hook length, determined by its relation to other cells of $\lambda$ and
$n(\lambda)$ is given by \[2n(\lambda)=\sum_{x\in\lambda}1+cn(x)-hl(x).\]

The coefficients of the series $\eval{\sum c_rt^r}=\sum\eval{c_r}t^r$ were
shown in chapter 4 of \cite{Aiston} to satisfy the simple recursive relation
\[\eval{c_{r+1}}=\frac{v^{-1}s^{-r}-vs^r}{s^{r+1}-s^{-r-1}}\eval{c_r},\]
leading to the product formula
\begin{equation}
\sum\eval{c_r}t^r=\prod_{i=0}^\infty \frac{
1+vs^{2i+1}t}{1+v^{-1}s^{2i+1}t}.
\label{Product}
\end{equation}

Applying Macdonald's formula (\ref{Schurmac}) gives 
\begin{equation}
\eval{\lambda} =s_\lambda{\eval{\sum
c_rt^r}}=q^{n(\lambda)}\prod_{x\in\lambda}\frac{vs-v^{-1}sq^{cn(x)}}{1-q^{hl(x)}}
\label{E-prodfrac}
\end{equation}
 taking $q=s^2$, $a=vs$ and $b=v^{-1}s$. This
can be rewritten in the form
\begin{equation}
\eval{\lambda}=\prod_{x\in\lambda}
\frac{v^{-1}s^{cn(x)}-vs^{-cn(x)}}{s^{hl(x)}-s^{-hl(x)}},\label{Schurprod}
\end{equation}
which can equally be derived by skein theory for
$\eval{Q_\lambda}$ from
\cite{AistonMorton}, without assuming that $Q_\lambda$ can be
expressed as a Schur function.

It is clear that the specialisation $v=s^{-N}$ gives
\begin{eqnarray}\eval{\lambda}_N=0&\iff&\mbox{some cell } x\in \lambda \mbox{
has content } -N \nonumber\\
&\iff&\lambda \mbox{ has }>N \mbox{ parts}. \label{LambdaSpec}
\end{eqnarray}

Consequently when $N\ge l(\lambda)$ we can always write
$t_\lambda(Y)_N=\displaystyle\frac{\hopf{\lambda}{Y}_N}{\eval{\lambda}_N}$.

\begin{lemma}
Let $\lambda$ be a Young diagram and let $N\geq l(\lambda)$.
Then the specialisation $E_{\lambda}^N(t)$ is a polynomial in $t$ of degree
$\le N$.
\end{lemma}
\begin{proof}
We must show that $t_\lambda(c_r)_N=0$ for $r>N$. Since ${\eval{\lambda}_N}\ne
0$ it is enough to show that $\hopf{\lambda}{c_r}_N=0$ for $r>N$.

By  (\ref{LambdaSpec}) we know that $\eval{c_r}_N=0$ for $r>N$, and so
\[\hopf{c_r}{\lambda}_N=t_{c_r}(\lambda)_N\eval{c_r}_N=0.\]
The result follows
since $\hopf{\lambda}{c_r}=\hopf{c_r}{\lambda}$.
\end{proof}

We now find a factorisation of this polynomial $E_\lambda^N(t)$ into linear
factors in lemma \ref{ElambdaN}, and use the Schur function
$s_\mu(E_\lambda^N(t))$ to derive a determinantal formula for
$\hopf{\lambda}{\mu}_N$ in theorem \ref{LambdaMu}.

Recall first the determinantal description of the Schur function 
\[s_\lambda(x_1,\ldots,x_N)=s_\lambda\left(\prod_{i=1}^N (1+x_it)\right).\]

For fixed $N\ge l(\lambda)$ define the index set \[I_\lambda=\{\lambda_1+N-1,
\lambda_2+N-2,\ldots,\lambda_N\}\] consisting of $N$ distinct integers $\ge0$.

Number the columns of the $N\times\infty$ matrix 
\[\left(\begin{array}{cccc}1&x_1&x_1^2&\cdots\\
1&x_2&x_2^2&\cdots\\\vdots&\vdots&\vdots&\vdots\\
1&x_N&x_N^2&\cdots
\end{array}\right)\] starting with $0$, so that its
 $(i,j)$ entry is 
$x_i^j$.

Write $P_\lambda^N({\bf x})$ for the $N\times N$ minor determined by the
columns $I_\lambda$.
 The leading minor is 
$P_\emptyset^N({\bf x})$, and the Schur function is given by
\begin{equation}
s_\lambda(x_1,\ldots,x_N)=P_\lambda^N({\bf
x})/P_\emptyset^N({\bf x}).
\label{Schurclass}
\end{equation}

The formula for $\hopf{\lambda}{\mu}_N$ can be found in terms of a function of
$q=s^2$ defined similarly by $N\times N$ minors of the infinite Vandermonde
matrix
\[V=\left(
      \begin{array}{cccccc}
         1 & 1 & 1 & 1 & 1 & \cdots\\
         1 & q &q^2&q^3&q^4& \cdots\\
         1 &q^2&q^4&q^6&q^8& \cdots\\
         1 &q^3&q^6&q^9&q^{12}&\cdots\\
         1 &q^4&q^8&q^{12}&q^{16}&\cdots\\
         \vdots&\vdots&\vdots&\vdots&\vdots&\ddots
      \end{array}
      \right)\]
 with $(i,j)$ entry $q^{ij}$, numbering both rows and columns from
$0$.  

\begin{definition}
For partitions $\lambda,\mu$ with
$l(\lambda), l(\mu)\ge N$  take index sets $I_\lambda$ and $I_\mu$ as above
and define $P^N_{(\lambda,\mu)}$ to be the $N\times N$ minor of $V$ with rows
indexed by $I_\mu$ and columns by $I_\lambda$.
\end{definition}

\begin{theorem} \label{LambdaMu}
\[\hopf{\lambda}{\mu}_N=s^{(1-N)(|\lambda|+|\mu|)}P^N_{(\lambda,\mu)}/
P^N_{(\emptyset,\emptyset)}.\]
\end{theorem}
\begin{proof}
Using  formula (\ref{Schurclass}) for the Schur function, applied to the
rows
$I_\mu$ of $V$, we can write
\begin{equation}
s_\lambda\left(\prod_{i\in
I_\mu}(1+q^it)\right)=P^N_{(\lambda,\mu)}/P^N_{(\emptyset,\mu)}.
\label{SchurLambdaMu}
\end{equation}

In particular
\[P^N_{(\lambda,\emptyset)}/P^N_{(\emptyset,\emptyset)}
=s_\lambda\left(\prod_{i=0}^{N-1}(1+q^it)\right).\]

We first derive a formula for the polynomial $E_\lambda^N(t)$ in the next
lemma.
\begin{lemma}\label{ElambdaN}
\[\prod_{i\in I_\lambda}(1+q^it)=E_\lambda^N(s^{N-1}t),\] giving
\[E_\lambda^N(t)=\prod_{j=1}^N(1+s^{N+2\lambda_j-2j+1}t).\]
\end{lemma}

To finish the proof of theorem \ref{LambdaMu} we use equation
(\ref{SchurLambdaMu}), lemma \ref{ElambdaN} and lemma \ref{L-muschur} to get
\begin{eqnarray*}
P^N_{(\lambda,\mu)}/P^N_{(\emptyset,\mu)}=s_\lambda\left(\prod_{i\in
I_\mu}(1+q^it)\right)&=&s_\lambda(E_\mu^N(s^{N-1}t))\\
&=&s^{(N-1)|\lambda|}s_\lambda(E_\mu^N(t))\\
&=& s^{(N-1)|\lambda|}\hopf{\lambda}{\mu}_N/\eval{\mu}_N
\end{eqnarray*}
and
\[P^N_{(\mu,\emptyset)}/P^N_{(\emptyset,\emptyset)}
=s_\mu(E_\emptyset^N(s^{N-1}t)
=s^{(N-1)|\mu|}\eval{\mu}_N.\]

Multiply the two expressions to get

\[P^N_{(\lambda,\mu)}/P^N_{(\emptyset,\emptyset)}=s^{(N-1)(|\lambda|+|\mu|)}
\hopf{\lambda}{\mu}_N.\]

The formula claimed in theorem \ref{LambdaMu} for $\hopf{\lambda}{\mu}_N$ in
terms of the
$N\times N$ minors of the Vandermonde matrix follows at once.
\end{proof}

\begin{proof}[of lemma \ref{ElambdaN}]

We begin by establishing theorem \ref{LambdaMu} directly in the case $\mu=c_k$.
  The index set $I_{c_k}$
is
\[I_{c_k}=\{N,N-1,\ldots,N-k+1,N-k-1,\ldots,2,1,0\}
=I_\emptyset\cup\{N\}-\{N-k\}.\]

Recall that  equation (\ref{E-prodfrac}) gives
\[
  E_{\emptyset}(t)=\sum_{r\geq 0}\eval{c_r}t^r=
  \prod_{i=0}^{\infty}\frac{1+vs^{2i+1}t}{1+v^{-1}s^{2i+1}t}.
\]
The substitution $v=s^{-N}$ reduces this to the finite product
\[  E^N_{\emptyset}(t)=\prod_{i=0}^{N-1}(1+s^{-N+2i+1}t).
\]
Then \begin{equation}
  E^N_{\emptyset}(s^{N-1}t)
                 =\prod_{i=0}^{N-1}(1+s^{2i}t)
=\prod_{i\in I_\emptyset}(1+q^it).
\end{equation}

Let $k$ be an integer, $k\leq N$. By corollary \ref{C-ckpower} we have
\[
  E_{c_k}(t)=\frac{1+v^{-1}st}{1+v^{-1}s^{-2k+1}t}E_{\emptyset}(t).
\]
Substituting $v=s^{-N}$ in the above equation gives
\begin{eqnarray}
  E^N_{c_k}(s^{N-1}t)&=&\frac{1+s^{2N}t}{1+s^{2N-2k}t}E^N_{\emptyset}(s^{N-1}t)
                \nonumber \\
              &=&\prod_{i\in I_{c_k}}(1+q^it).
\end{eqnarray}

Now  $s_\lambda\left(\prod_{i\in
I_{c_k}}(1+q^it)\right)=P_{(\lambda,c_k)}^N/P_{(\emptyset,c_k)}^N$, by
(\ref{SchurLambdaMu}).
Consequently
\[P_{(\lambda,c_k)}^N/P_{(\emptyset,c_k)}^N=s_\lambda(E^N_{c_k}(s^{N-1}t))\]
and
\[P_{(c_k,\emptyset)}^N/P_{(\emptyset,\emptyset)}^N=
s_{c_k}(E^N_{\emptyset}(s^{N-1}t)).\]
Since
$s_\lambda(E^N_{c_k}(t))=\hopf{\lambda}{c_k}_N/\eval{c_k}_N$ by lemma
\ref{L-muschur}, we can multiply the two expressions to get
\begin{eqnarray}
P_{(\lambda,c_k)}^N/P_{(\emptyset,\emptyset)}^N
&=&s_\lambda(E^N_{c_k}(s^{N-1}t)) s_{c_k}(E^N_{\emptyset}(s^{N-1}t))\nonumber\\
&=&s^{(N-1)|\lambda|}\frac{\hopf{\lambda}{c_k}_N}{\eval{c_k}_N}
s^{(N-1)|c_k|}\eval{c_k}_N \nonumber\\
&=&s^{(N-1)(|\lambda|+|c_k|)}\hopf{\lambda}{c_k}_N.\label{Lambdack}
\end{eqnarray} 

To complete the proof of lemma \ref{ElambdaN} it is enough to show that
\[s_{c_k}\left(\prod_{i\in
I_\lambda}(1+q^it)\right)=s_{c_k}(E^N_{\lambda}(s^{N-1}t))\] for all $k\le N$,
so as to compare the coefficients of the two polynomials.

By 
(\ref{SchurLambdaMu}) we have $s_{c_k}\left(\prod_{i\in
I_\lambda}(1+q^it)\right)=P_{(c_k,\lambda)}^N/P_{(\emptyset,\lambda)}^N$, while
\begin{eqnarray*}
s_{c_k}(E^N_{\lambda}(t))&=&\hopf{\lambda}{c_k}_N/\eval{\lambda}_N\\
&=&s^{(1-N)|c_k|}P_{(\lambda,c_k)}^N/P_{(\lambda,\emptyset)}^N
\end{eqnarray*} by
(\ref{Lambdack}).

Since $P_{(c_k,\lambda)}^N=P_{(\lambda,c_k)}^N$ and $
s_{c_k}(E^N_{\lambda}(s^{N-1}t))=s^{(N-1)|c_k|}s_{c_k}(E^N_{\lambda}(t))$ the
result follows.
\end{proof}

 We now deduce a
formula for the power series
$E_{\lambda}(t)$ with coefficients in the two variables $s$ and $v$ from the
formula for
$E_{\lambda}^N(t)$ in lemma \ref{ElambdaN}.

\begin{theorem}\mylabel{T-elambdaXprod}
We have
\[
  E_{\lambda}(t)=\prod_{j=1}^{l(\lambda)}
   \frac{1+v^{-1}s^{2\lambda_j-2j+1}t}
        {1+v^{-1}s^{-2j+1}t}E_{\emptyset}(t)
\]
for any Young diagram $\lambda$.
\end{theorem}
\begin{proof}
For any integer $N\geq l(\lambda)$ we have
$E_{\lambda}^N(s^{N-1}t)=\prod_{i\in I_\lambda}(1+q^it)$ by lemma
\ref{ElambdaN}.
Now
\[I_\lambda=I_\emptyset\cup_{j=1}^{l(\lambda)}\{\lambda_j+N-j\}-
\cup_{j=1}^{l(\lambda)}\{N-j\}.\]

Hence
\[
E_{\lambda}^N(s^{N-1}t)
=
\prod_{j=1}^{l(\lambda)}\frac{1+q^{\lambda_j+N-j}t}{1+q^{N-j}t}
E_{\emptyset}^N(s^{N-1}t)
\] and so
\begin{eqnarray*}
E_{\lambda}^N(t)
&=&
\prod_{j=1}^{l(\lambda)}\frac{1+s^{2(\lambda_j+N-j)}s^{1-N}t}
{1+s^{2(N-j)}s^{1-N}t}E_{\emptyset}^N(t)\\
&=&
  \prod_{j=1}^{l(\lambda)}\frac{1+s^{N+2\lambda_j-2j+1}t}
                               {1+s^{N-2j+1}t} E_{\emptyset}^N(t).
\end{eqnarray*}

This means that the two power series $E_{\lambda}(t)$
and
\[
  \prod_{j=1}^{l(\lambda)}
  \frac{1+v^{-1}s^{2\lambda_j-2j+1}t}
        {1+v^{-1}s^{-2j+1}t}E_{\emptyset}(t)
\]
agree for every substitution $v=s^{-N}$ with
$N\geq l(\lambda)$.
Since the coefficients are Laurent polynomials in $v$ it follows that the two
power series are equal when $v$ is treated as an indeterminate.
\end{proof}

When we apply theorem \ref{T-elambdaXprod} to the case $\lambda=c_k$
and compare the result with corollary \ref{C-ckpower}
we note a number of cancellations in 
\[
  \prod_{j=1}^{l(\lambda)}
  \frac{1+v^{-1}s^{2\lambda_j-2j+1}t}
       {1+v^{-1}s^{-2j+1}t}.
\]
In our final theorem we give a simpler expression for the rational
function $E_\lambda(t)/E_\emptyset(t)$, in terms of the Frobenius notation
$\lambda=(a_1,\ldots,a_{d(\lambda)}|b_1,\ldots,b_{d(\lambda)})$ for the
partition
$\lambda$. As in
\cite{Mac} recall that the Frobenius notation describes $\lambda$ in terms of
the lengths $\{a_i\}$ of the `arms' and $\{b_i\}$ of the `legs' of its Young
diagram, counting right or down from the $d(\lambda)$ cells on the main
diagonal. Thus
$a_i=\lambda_i-i$ 
 and
$b_i=\trans{\lambda}_i-i$ for $i=1,\ldots,d(\lambda)$.

\begin{theorem}\label{Frobenius}
We have
\[
  E_{\lambda}(t)=
  \prod_{i=1}^{d(\lambda)}\frac{1+v^{-1}s^{2a_i+1}t}
               {1+v^{-1}s^{-2b_i-1}t}E_\emptyset(t)
\]
for any Young diagram
$\lambda=(a_1,\ldots,a_{d(\lambda)}|b_1,\ldots,b_{d(\lambda)})$ in Frobenius
notation.
\end{theorem}
\begin{proof}
The result follows from exactly the same argument as in theorem
\ref{T-elambdaXprod} once we observe that we can write the index set
$I_\lambda$ as
\[I_\lambda=I_\emptyset\cup_{i=1}^{d(\lambda)}\{N+a_i\}
-\cup_{j=1}^{d(\lambda)}\{N-b_j-1\}.\]

To confirm this, write $a_i=\lambda_i-i, 
b_j=\trans{\lambda}_j-j$ for all $i,j\le N$. Then
\begin{eqnarray*} a_i\ge 0&\iff&i\le d_\lambda,\\
b_j\ge0&\iff& j\le d_\lambda.
\end{eqnarray*}
The integers $a_i, 1\le i\le N$ are all distinct, and $I_\lambda=\{N+a_i, i\le
N\}$. Since $N+a_i\ge N\iff i\le d_\lambda$ there are exactly $d_\lambda$
numbers in $I_\lambda$ which are not in $I_\emptyset$, and hence $d_\lambda$
 numbers in $I_\emptyset$ which are not in  $I_\lambda$. We can identify
these as $\{N-b_j-1:j\le d_\lambda\}$ by showing that
\[N-b_j-1\ne N+a_i, \mbox{ for any } i,j\le N.\]
Consider the cell in row $i$ and column $j$, $ 1\le i,j\le N$, in relation to
the Young diagram $\lambda$. If the cell lies in $\lambda$ then $\lambda_i\ge
j$ and $\trans{\lambda}_j\ge i$. Add these to get $a_i+b_j\ge 0$ and so
$a_i+b_j\ne -1$. Otherwise  $\lambda_i\le
j-1$ and $\trans{\lambda}_j\le i-1$, giving  $a_i+b_j\le -2$ and
again
$a_i+b_j\ne -1$.

This completes the proof;  indeed  the rational function
\[\prod_{i=1}^{d(\lambda)}\frac{1+v^{-1}s^{2a_i+1}t}
               {1+v^{-1}s^{-2b_i-1}t}=E_{\lambda}(t)/E_\emptyset(t)\]
admits no further cancellation, since the exponents of $s$ in the numerator
are all positive, while those in the denominator are negative.
\end{proof}

Having used this result to get the power series for $E_\lambda(t)$ with
explicit coefficients in $\Lambda$, the general invariant $\hopf{\lambda}{\mu}$
is then found by calculating  Schur functions of the series, using the
formula \[\hopf{\lambda}{\mu}=s_\mu(E_\lambda(t))\eval{\lambda}.\]

\begin{remark} The rational function above can
be expressed in terms of the {\em content polynomial} $C_\lambda(t)$ defined
by
\[C_\lambda(t)=\prod_{x\in\lambda}(1+q^{cn(x)}t).\] It is easy to check,
by decomposing $\lambda$ as a disjoint union of simple hooks based on the
diagonal cells, that \[\prod_{i=1}^{d(\lambda)}\frac{1+v^{-1}s^{2a_i+1}t}       
{1+v^{-1}s^{-2b_i-1}t}=C_{\lambda}(v^{-1}st)/C_{\lambda}(v^{-1}s^{-1}t).\]
 This leads to an alternative proof of theorem \ref{Frobenius} using 
theorem 3.10 of \cite{Murphy}, once it is established that the product of
the Gyoja-Aiston idempotent
$e_\lambda$, with $|\lambda|=n$, and the polynomial $EM(t)$, whose
coefficients are the elementary symmetric functions of the Murphy operators
in $H_n$,  satisfies the equation
\[EM(t)e_\lambda=C_\lambda(t)e_\lambda.\]
\end{remark}

\subsection{An explicit example} \label{JacobiTrudy}
We illustrate our results by calculating the framed Homfly polynomial
$\hopf{\Threeone}{\Twotwo}$.

 Let $\mu$ be a partition whose conjugate  $\trans{\mu}$
has
$r$ parts. The Jacobi-Trudy formula  gives the Schur function $s_\mu$ of the
series $\sum e_it^i$  as the determinant of the $r\times r$ matrix
whose
$(i,j)$ entry is
$e_k$,  with $k=\trans{\mu}_i+j-i$, taking $e_k=0$ for $k<0$. The
diagonal entries in the matrix are thus the coefficients corresponding to
the columns of $\mu$, and each row of the matrix is completed by taking
consecutive coefficients from the series. Then
\[s_\Twotwo=\det\left(\begin{array}{cc}e_2&e_3\\e_1&e_2
\end{array}\right)=e_2^2-e_1e_3\] and it is enough to expand the series
$E_\Threeone(t)$ as far as the term in $t^3$ in order to calculate
$s_\Twotwo(E_\Threeone(t))$.

For the diagram $\lambda=\Threeone$ we have $d(\lambda)=1$ with $a_1=2,b_1=1$
so that \[ E_{\Threeone}(vs^{-1}t)=
  \frac{1+q^2t}
               {1+q^{-2}t}E_\emptyset(vs^{-1}t).\]
Using (\ref{Schurprod}) with $\lambda=c_i$ we can write
\begin{eqnarray*}E_\emptyset(vs^{-1}t)&=&1+\frac{1-v^2}{q-1}t+\frac{(1-v^2)
(1-qv^2)}{(q-1)(q^2-1)}t^2\\
&+&\frac{(1-v^2)(1-qv^2
)(1-q^2v^2)}{(q-1)(q^2-1)(q^3-1)}t^3+O(t^4).\end{eqnarray*}

We have
\[\frac{1+q^2t}{1+q^{-2}t}=1+(q^2-q^{-2})t-(1-q^{-4})t^2+(q^{-2}-q^{-6})t^3+O(t^4)\]
and so
\begin{eqnarray*}E_{\Threeone}(vs^{-1}t)&=&1+\left(\frac{1-v^2}{q-1}+q^2-q^{-2}\right)t
\\
&+&\left(\frac{(1-v^2)(1-qv^2)}{(q-1)(q^2-1)}+
(q^2-q^{-2})\frac{1-v^2}{q-1}-(1-q^{-4})\right)t^2\\
&+&\left(\frac{(1-v^2)(1-qv^2)(1-q^2
v^2)}{(q-1)(q^2-1)(q^3-1)}+(q^2-q^{-2})\frac{(1-v^2)(1-qv^2)}{(q-1)(q^2-1
)}\right.\\
&-&\left. (1-q^{-4})\frac{1-v^2}{q-1}+(q^{-2}-q^{-6})\right)t^3+O(t^4)
  \\
&=& 1+e_1t+e_2t^2+e_3t^3+O(t^4).
\end{eqnarray*}

Now
\begin{eqnarray*}s_\Twotwo(E_\Threeone(t))&=&
(vs^{-1})^{-4}s_\Twotwo(E_\Threeone(vs^{-1}t))\\
&=&(vs^{-1})^{-4}(e^2_2-e_1e_3),
\end{eqnarray*} where $e_1, e_2$ and $e_3$ are the coefficients in the series
above.

Combined with the expression 
\[\eval{\Threeone}=
\frac{(v^{-1}-v)(v^{-1}s-vs^{-1})(v^{-1}s^2-vs^{-2})(v^{-1}s^{-1}-vs)}
{(s-s^{-1})^2(s^2-s^{-2})(s^4-s^{-4})}\]
from 
(\ref{Schurprod}) we get
\begin{eqnarray}
\hopf{\Threeone}{\Twotwo}&=&
(vs^{-1})^{-4}(e^2_2-e_1e_3)\eval{\Threeone}\nonumber\\
&=&
v^{-8}\frac{(v^2-1)^2(v^2-q)(v^2-q^2)(v^2q-1)}{(q-1)^3(q^2-1)^3(q^3-1)(q^4-1)}C
\label{Hopfex}\end{eqnarray}
with
\begin{eqnarray*}C&=&
 - q^{14} + q^{13} + q^{12} - q^{11} - q^{10} - q^{9} + 2q^{8} + q^{7} -
q^{6}  - 2q^{4} + 2q^{2} - 1\\
&& + (q^{13} - q^{10} - q^{9} + 2q^{8} + q^{7}+ q^{6} - 2q^{4} + q^{2} +
q)v^{2}\\ &&+ ( - q^{10} - q^{9} + q^{8} - q^{4} -q^{3})v^{4}+ q^{6}v^{6}.
\end{eqnarray*}

\begin{remark}
Since  $\hopf{\lambda}{\mu}=s_\mu(E_\lambda(t))\eval{\lambda}
=s_\lambda(E_\mu(t))\eval{\mu}$ the calculation in our example could have
been done with the roles of $\lambda$ and $\mu$ interchanged. In this case we
would have had a more complicated expression for $E_\mu(t)$ since
$d(\mu)=2$, and in addition the expression for $s_\lambda$ would involve terms
up to degree 4 from the Jacobi-Trudy formula
\[s_\Threeone=\det\left(\begin{array}{ccc}e_2&e_3&e_4\\1&e_1&e_2\\0&1&e_1
\end{array}\right).\]

 In general $\hopf{\lambda}{\mu}$  contains some form of least
common multiple of
$\eval{\lambda}$ and
$\eval{\mu}$, and a `core' part, such as the expression $C$ in the example
above. 
\end{remark}

To find $\hopf{\lambda}{\mu}_N$ we can make the substitution $v=s^{-N}$. When
either $\eval{\lambda}_N=0$ or $\eval{\mu}_N=0$ this will vanish. This happens
when one of the partitions has more than $N$ parts. In our example this is the
case when $N=1$, as we can see from the factor $v^2q-1$ in (\ref{Hopfex}).

We can  equally calculate $\hopf{\lambda}{\mu}_N$
from  Vandermonde minors, using theorem
\ref{LambdaMu}, and we compare the two results for our example in the case
$N=3$.

The substitution with $N=3$ in equation (\ref{Hopfex}) gives
\[\hopf{\Threeone}{\Twotwo}_3=(q^2+q+1)(q^8+q^4+q^3-q^2+1)(q^2+1)(q^4+q^
3+q^2+q+1)/q^3.\]

On the other hand  taking
$I_\Threeone=\{0,2,5\}$ and
$I_\Twotwo=\{0,3,4\}$, gives the minor
\begin{eqnarray*}P^3_{(\Threeone,\Twotwo)}&=&\det\left(\begin{array}{ccc}
1&1&1\\1&q^6&q^{15}\\1&q^8&q^{20}
\end{array}\right)=q^{26}-q^{23}-q^{20}+q^{15}+q^8-q^6\\
&=&q^{6}(q + 1)(q^{2} + 1)(q^{2} + q +
1)\\
&&\quad (q^{4} + q^{3} + q ^{2} + q + 1)(q^{8} +
q^{4} + q^{3} - q^{2} + 1)(q -
1)^{3},
\end{eqnarray*} and so
\begin{eqnarray*}\hopf{\Threeone}{\Twotwo}_3&=&
q^{-8}P^3_{(\Threeone,\Twotwo)}/P^3_{(\emptyset,\emptyset)}, \mbox{ by 
 theorem \ref{LambdaMu},}\\ &=&
q^{-8}\frac{q^{26}-q^{23}-q^{20}+q^{15}+q^8-q^6}
{(q-1)(q^2-1)(q^2-q)}\\
&=&(q^2+1)(q^2+q+1)\\
&&\quad (q^4+q^
3+q^2+q+1)(q^8+q^4+q^3-q^2+1)/q^3, 
\end{eqnarray*}
as before.
\subsection{The $\mysln_q$ invariants of the Hopf link}
Irreducible $\mysln_q$ modules$\{V_\lambda\}$ are indexed by Young diagrams
with at most $N$ rows. The framed $\mysln_q$ invariant $(J_H;V_\lambda,V_\mu)$
for the Hopf link $H$ coloured by the irreducible modules $V_\lambda$ and
$V_\mu$ is given, up to a fractional power of $s$, from our Homfly invariant
$\hopf{\lambda}{\mu}$ by applying the $\mysln$ specialisation $v=s^{-N}$. The
exact formula, noted in \cite{Lukac},  is
\[(J_H;V_\lambda,V_\mu)=x^{2|\lambda||\mu|}\hopf{\lambda}{\mu}_N, \mbox{ with
}x=s^{-\frac{1}{N}}.\]
The correction factor arises from the isomorphism between the Hecke algebras
used here and the version based on the $R$-matrix in $\mysln_q$. 

The $\mysln_q$ invariants of the Hopf link are then essentially the Laurent
polynomials in $q$ determined by $N\times N$ minors of the Vandermonde matrix
$(q^{ij})$.

In the simplest case $N=2$ the most general minor is
\[\det\left(\begin{array}{cc}
q^{ij}&q^{i(j+b)}\\
q^{(i+a)j}&q^{(i+a)(j+b)}
\end{array}\right)=q^{2ij+ib+aj}(q^{ab}-1),\] arising from partitions
$\lambda$ with parts $(b+j-1,j)$ and $\mu$ with parts $(a+i-1,i)$.

Then $\displaystyle \hopf{\lambda}{\mu}_2=\frac{q^{ab}-1}{q-1}$, up to a power
of
$q$. Calculation of this power and the correction factor gives the formula
\begin{equation}(J_H;V_\lambda,V_\mu)=[ab]\label{SLtwo}\end{equation}
for the $sl(2)_q$ invariants of the Hopf link,  where
$[k]$ is the quantum integer $\displaystyle [k]=\frac{s^k-s^{-k}}{s-s^{-1}}$.

The $sl(2)_q$ modules $V_\lambda$ and $V_\mu$ in this case do not depend on
$i$ or $j$, and have dimension $a$ and $b$ respectively. They correspond to
the single row diagrams $d_{a-1}$ and $d_{b-1}$. The
formula (\ref{SLtwo}) for the $sl(2)_q$ Hopf link invariants was given in
\cite{MorStr} and
\cite{KM}. It allows a simple development of the properties for the related
3-manifold invariants where $q$ is replaced by a root of unity, \cite{MorStr2,
KM}.

For $N\ge 3$ the $\mysln_q$ invariants of the Hopf link do not have such an
immediately memorable form. Their expression by Vandermonde minors was 
given by Kohno and Takata, \cite{Kohno}, only in the case where $q$ is a root
of unity. Our determination of them here for generic $q$ in theorem
\ref{LambdaMu} and the eventual 2-variable formulation for
$\hopf{\lambda}{\mu}$ given in theorem
\ref{Frobenius} was initially inspired by our reading of \cite{Kohno}.

\medskip
\noindent
Version 1.6, July 2001.

Copies of items marked * can be found on the Liverpool Knot Theory site
{\tt http://www.liv.ac.uk/\~{}su14/knotgroup.html} or by following links from
\break
{\tt http://www.liv.ac.uk/maths/} .
\end{document}